%% file: main.tex
\documentclass[letterpaper, 10 pt, conference]{ieeeconf}   
\IEEEoverridecommandlockouts  %
\usepackage{amsmath}

\usepackage{amsthm}
\usepackage{amssymb}  
\usepackage{mathtools}
\usepackage{bm}

\usepackage{url}
\usepackage{import}

\usepackage[dvipsnames]{xcolor}
\usepackage{tabularx} 
\usepackage{multirow}
\usepackage{cite}
\usepackage{algorithm} 
\usepackage{algpseudocode} 
\usepackage{empheq}
\usepackage{multirow}
\usepackage{float}

\usepackage{booktabs}
\usepackage{tikz}

\makeatletter
\let\MYcaption\@makecaption
\makeatother

\usepackage{subcaption}
\usepackage{overpic}

\makeatletter
\let\@makecaption\MYcaption
\makeatother

\usepackage{balance} 

\let\labelindent\relax %
\usepackage{enumitem}

\makeatletter
\let\NAT@parse\undefined
\makeatother
\usepackage{hyperref} 
\hypersetup{
    colorlinks=true,
    citecolor=black,
    linkcolor=black,
    urlcolor=black,
    hypertexnames=true
}

\newtheorem{thm}{Theorem}

 \theoremstyle{definition}
\newtheorem{defn}[thm]{Definition}

\newtheorem{rem}[thm]{Remark}

\let\pf\proof
\let\endpf\endproof
\let\citep\cite

\allowdisplaybreaks

\input{macros}

\newcommand{\IEEEacceptance}{
    \begin{tikzpicture}[overlay, remember picture]
        \path (current page.north east) ++(-1.95,-0.2) node[below left] {
            This paper has been accepted for publication in the proceedings of the IEEE 63rd Conference on Decision and Control.
        };
    \end{tikzpicture}
}

\newcommand\copyrighttext{%
  \footnotesize \textcopyright 2024 IEEE. Personal use of this material is permitted.
  Permission from IEEE must be obtained for all other uses, in any current or future 
  media, including reprinting/republishing this material for advertising or promotional 
  purposes, creating new collective works, for resale or redistribution to servers or 
  lists, or reuse of any copyrighted component of this work in other works.}
\newcommand\copyrightnotice{%
\begin{tikzpicture}[remember picture,overlay]
\node[anchor=south,yshift=10pt] at (current page.south) {\fbox{\parbox{\dimexpr\textwidth-\fboxsep-\fboxrule\relax}{\copyrighttext}}};
\end{tikzpicture}%
}

\title{\LARGE \bf
    Compositional Construction of Barrier Functions for Switched Impulsive Systems
}

\author{%
    Katharina Bieker$^{1,*}$, %
    Hugo Kussaba$^{2}$, %
    Philipp Scholl$^{1}$, %
    Jaesug Jung$^{2}$, \\ %
    Abdalla Swikir$^{2}$, %
    Sami Haddadin$^{2}$, %
    and %
    Gitta Kutyniok$^{1}$
\thanks{KB, PS, JJ, SH, and GK were supported by LMUexcellent and TUM AGENDA 2030, funded by the Federal Ministry of Education and Research (BMBF) and the Free State of Bavaria under the Excellence Strategy of the Federal Government and the Länder as well as by the Hightech Agenda Bavaria. GK was also supported in part by the DAAD programme Konrad Zuse Schools of Excellence in Artificial Intelligence, sponsored by the Federal Ministry of Education and Research. 
HK received support from the European Union's Horizon 2020 research and innovation programme (Marie Skłodowska-Curie grant agreement no. 899987).
}%
\thanks{$^{1}$ Department of Mathematics at Ludwig-Maximilians-Universität München, Germany. G. Kutyniok is also with the Munich Center for Machine Learning (MCML), Germany.}%
\thanks{$^{2}$ Munich Institute of Robotics and Machine Intelligence (MIRMI), Technical University of Munich (TUM), Germany. A. Swikir is also with the Department of Electrical and Electronic Engineering, Omar Al-Mukhtar University (OMU), Albaida, Libya.}%
\thanks{$^{*}$ Corresponding author:
{\tt\small bieker@math.lmu.de}}%
}

\begin{document}
	
\maketitle
\IEEEacceptance
\copyrightnotice
\thispagestyle{plain}
\pagestyle{plain}

\begin{abstract}
    Many systems occurring in real-world applications, such as controlling the motions of robots or modeling the spread of diseases, are switched impulsive systems. To ensure that the system state stays in a safe region (e.g., to avoid collisions with obstacles), barrier functions are widely utilized. As the system dimension increases, deriving suitable barrier functions becomes extremely complex.
    Fortunately, many systems consist of multiple subsystems, such as different areas where the disease occurs. In this work, we present sufficient conditions for interconnected switched impulsive systems to maintain safety by constructing local barrier functions for the individual subsystems instead of a global one, allowing for much easier and more efficient derivation.
    To validate our results, we numerically demonstrate its effectiveness using an epidemiological model.
\end{abstract}

\input{introduction}

\input{review}
\input{results}

\input{simulation}

\section{Conclusion}\label{sec:Conc}
In this paper, we address the challenge of constructing barrier functions in practice, especially for high-dimensional systems, by considering a compositional construction for interconnected systems. 
By exploiting fundamental principles, we present a novel formulation for the (compositional) construction of barrier functions in interconnected switched impulsive systems. Thereby, we could avoid inferring a small-gain condition and were able to derive conditions that allow for a parallelized computation of the local pseudo barrier functions, cf.\ Theorem~\ref{th:Adapted_Inequality}. To illustrate our approach, we provide a numerical example based on the SIR model, which proves that the compositional computation saves a considerable amount of computing time. Moreover, it becomes clear that in this specific example, a global construction of the barrier function is not possible with means of the SOS-method as the unsafe set $\unsafe$ can not be characterized by a polynomial, whereas the problem can be solved using the compositional construction. 
An important extension of our formulation would be to utilize time-varying barrier functions as done in \cite{TN2023}. This way, control laws depending on the time (and not only on the state) can be considered. Moreover, the time between jumps and impulses can be taken into account, allowing for more generally composed barrier functions. %

\footnotetext{The experiments were performed using the optimization solver MOSEK 10.1 \citep{MOSEK} and on a notebook equipped with an AMD Ryzen 7 5700U CPU and 16 GB of RAM. %
}

\newpage

\balance

\end{document}

%% file: macros.tex
\newcommand{\N}{\mathbb{N}}
\newcommand{\R}{\mathbb{R}}

\newcommand{\safe}{X_{\mathsf{safe}}}
\newcommand{\unsafe}{X_{\mathsf{unsafe}}}

\newcommand{\BsumP}[0]{B_{\mathsf{sum}}}
\newcommand{\BmaxP}[0]{B_{\mathsf{max}}}
\newcommand{\Bsum}[1]{B_{\mathsf{sum}}\left(#1\right)}
\newcommand{\Bmax}[1]{B_{\mathsf{max}}\left(#1\right)}
\newcommand{\faa}{\text{f.a.a. }} %

%% file: introduction.tex
\section{Introduction}  

Dynamical systems need to be controlled in a wide range of applications. These systems are often safety-critical, for instance, in aviation, healthcare, and industrial processes. 
To ensure that such control systems are safe, i.e., stay within given boundaries, safety verification techniques such as reachability analysis have been proposed to prove that, starting from specific initial conditions, a system state cannot navigate into unsafe regions.
To avoid the explicit computation of the system trajectories, so-called \emph{barrier functions} are employed, an important and widely used tool for proving safety in dynamical systems \citep{prajna2007framework}. %

However, computing barrier functions in practical applications is usually difficult, especially for high-dimensional systems, as the computational costs grow with increasing dimension. 
The authors of \cite{jagtap2018temporal} and \cite{wongpiromsarn2015automata} were able to overcome the curse-of-dimensionality by searching for a polynomially parameterized barrier function but still suffered from a polynomially growing computing time with respect to the dimension of the system.
Thus, general approaches, such as those presented in \cite{jagtap2020formal}, become computationally intractable for large-scale systems.
Fortunately, many real-world systems can be modelled as \emph{interconnected systems}, i.e., as a collection of multiple smaller subsystems that are connected to each other, allowing for efficient compositional construction of barrier functions by so-called local barrier functions \citep{jagtap2020compositional,LYU22022, NSZ2022}.

In this work, we aim at the (compositional) construction of barrier functions to allow for the computational verification of whether a controlled \emph{switched impulsive systems} is safe. The detailed setting, along with the introduction of barrier functions, is presented in Section~\ref{sec:Basics}.
Such systems belong to the class of so-called \emph{hybrid systems} \citep{Lin2022,guan2005hybrid} and occur in many different practical applications, for instance, in robotics \citep{ding2011hybrid,ossadnik2022bsa}, biological systems as disease modeling \citep{wang2015stochastic,gao2018analysis}, or manufacturing processes \citep{CPW2001}. 

There already exist different versions of barrier functions for switched systems as well as for impulsive systems and hybrid systems \cite{KHS+2013, KKW2019, FBM2020}.
Building upon these theoretical foundations, we derive conditions for the compositional construction of a barrier function based on local pseudo-barrier functions in Section~\ref{sec:Main}. To the best of our knowledge, this is the first formulation for the composed construction of barrier functions for interconnected switched impulsive systems.
Our ideas for the construction are essentially based on the results from \cite{jagtap2020compositional}, \cite{NSZ2022}, and \cite{dashkovskiy2011small}.
However, in contrast to past works (e.g. \citep{NSZ2020, NSZ2022}), our derived conditions enable us to compute the local barrier functions of the subsystems in parallel, as they do not depend on any outside information from the other subsystems. 
To verify the derived results and demonstrate the reduction of computing time through the composite construction of the barrier functions, in Section~\ref{sec:Numerics}, we consider as a numerical example the SIR model \cite{liu2017infectious}, an epidemiological model that, for instance, was used to predict the infections during the COVID-19 pandemic.

%% file: review.tex
\section{Theoretical Basics}\label{sec:Basics}

\subsection{Barrier Functions}\label{subsec:BarrierFunctions}

    To formalize the concept of barrier functions, we first introduce the exact setting of the switched impulsive systems that we consider in this work. 

    \begin{defn}
    \label{def:system}
        A \emph{switched impulsive system} is given by a tuple 
        \begin{align}
            \Sigma=(X,U,\Omega, (f_p)_{p\in \{1,\dots,N\}}, g),
        \end{align}
        where 
        \begin{itemize}
            \item the \emph{state space} $X \subseteq \R^n$ is open and connected,  
            \item the \emph{control set} $U \subseteq \R^{n_u}$ is a compact set,
            \item the \emph{transition maps} $f_p\colon X \times U \rightarrow \R^n$, $p \in \{1,\dots,N\}$, are $C^1$ in the first and $C^0$ in the second argument, and 
            \item $\Omega=\{0<t_1<t_2<\dots\}$ is the set of points in time at which the system is subject to an impulse given by the \emph{jump map} $g\colon X \times U \times \R_{\ge0} \rightarrow X$. Assume that there exists $\delta > 0$ such that $t_i -t_{i-1} > \delta$ for all $i \in \N$ to ensure that there are only finite many jumps in a fixed time interval, avoiding Zeno solutions \citep{Dashkovskiy2017}.  
            
        \end{itemize}

        Furthermore, the following system equations have to hold:
        \begin{equation}\label{eq:system_control}
        \begin{aligned}
            \Sigma:
            \left\{
            \begin{array}{lll}
                \dot{\mathbf{x}}(t) &= f_{\boldsymbol{\sigma}(t)}(\mathbf{x}(t),\mathbf{u}(t)), \quad &\text{if } t \in \R_{\ge 0} \setminus \Omega,\\
                \mathbf{x}(t) &= g(\mathbf{x}^-(t),\mathbf{u}^-(t),t), \quad &\text{if } t\in \Omega,
            \end{array}
            \right.
        \end{aligned}
        \end{equation}
        where 
        \begin{itemize}
            \item ${\mathbf{x}}: [0,\infty) \to X$ is the \textit{state trajectory},
            \item the \emph{control input function} $\mathbf{u}\colon [0,\infty) \rightarrow U$ is measurable and right-continuous with left limits for all $t \in \R_{\ge0}$, and
            \item $\boldsymbol{\sigma} \colon [0,\infty) \rightarrow\{1,\dots,N\}$ is an admissible \emph{switching signal}, i.e., a piecewise constant function, right-continuous and the number of instants where the signal switches is finite on every finite interval (no Zeno behavior).
        \end{itemize}
        Moreover, by $(\cdot)^-$, we denote the limit from the left, i.e., 
        \[\mathbf{x}^-(t) = \displaystyle \lim_{s\nearrow t}\mathbf{x}(s) \quad \text{and} \quad \mathbf{u}^-(t) = \displaystyle \lim_{s\nearrow t}\mathbf{u}(s).\]
       We assume that the sets $\safe, \unsafe \subseteq X$ are a partition of $X$, i.e., $\safe \cap \unsafe = \emptyset$ and ${\safe \cup \unsafe = X}$  (defining the regions where the system state is safe or unsafe depending on the application).
        By $X_0 \subseteq \safe$, we denote the set of initial states. 
        In addition, we assume that a solution ${\varphi_{x_0,\mathbf{u}, \boldsymbol{\sigma}}(t)\colon [0,T)\rightarrow X}$ of \eqref{eq:system_control} that is absolutely continuous for all $t \in \R_{\ge 0}\setminus \Omega$ exists for any initial value $x_0 \in X_0$ and all $T >0$.%
    \end{defn}

    Following, we introduce sufficient conditions for barrier functions of switched impulsive systems to ensure that the controlled system is safe. The theorem is composed of the results presented in \cite{KHS+2013}, \cite{KKW2019}, and \cite{FBM2020}. 
     
    \begin{thm}\label{th:barrier_basic}
        Consider a switched impulsive system $\Sigma$ as in Definition~\ref{def:system}. Let $\mathbf{u}\colon X \rightarrow U$ be a continuous control law and $\boldsymbol{\sigma}\colon [0,\infty) \rightarrow \{1,\dots,N\}$ a predefined switching signal. Furthermore, assume that there exists a barrier function $B\colon X \rightarrow \R$ that is $C^1$ (almost everywhere) and satisfies the following conditions:
        \begin{subequations}
        \label{eq:Barrier}
        \begin{alignat}{2}
            &B(x) \le 0, && \forall x \in X_0,  
            \label{eq:Barrier_safe} \\
            &B(x) > 0, &&\forall x \in \unsafe, 
            \label{eq:Barrier_unsafe}\\
            &\frac{\partial B}{\partial x} f_p(x,\mathbf{u}(x)) \le \lambda B(x), \quad && \faa x \in X, \forall p \in \{1,\dots,N\},
            \label{eq:Barrier_dot}\\ 
            &B(g(x,\mathbf{u}(x),t)) \le \theta B(x), &&\forall x \in \overline{\safe}, t \in \Omega,
            \label{eq:Barrier_impulse}
        \end{alignat}
        \end{subequations}
        for some $\lambda \in \R$ and $\theta \in \R_{\ge0}$, where $\overline{\safe}$ denotes the closure of the set $\safe$.
        Then, for every $T>0$, the system $\Sigma$ is safe with respect to the control function $\mathbf{u}$, i.e., ${\varphi_{x_0,\mathbf{u}, \boldsymbol{\sigma}}(t) \in \safe}$ for all $t \in [0,T)$ for every solution $\varphi_{x_0,\mathbf{u}, \boldsymbol{\sigma}}\colon [0,T)\rightarrow X$ of \eqref{eq:system_control} with initial value ${\varphi_{x_0,\mathbf{u}, \boldsymbol{\sigma}}(0) = x_0 \in X_0}$ (for the control law $\mathbf{u}\colon X\rightarrow U$ and the switching signal $\boldsymbol{\sigma}\colon [0,\infty)\rightarrow \{1,\dots,N\}$).
    \end{thm}    
    \begin{pf}
        First, we define the function $h_p\colon X \rightarrow \R$ for all ${p \in \{1,\dots,N\}}$ by
        \begin{align*}
            h_p(x) := 
            \begin{cases}
                \frac{\partial B}{\partial x}(x) f_p(x,\mathbf{u}(x)) - \lambda B(x), & \text{if}~\frac{\partial B}{\partial x}(x)~\text{exists},\\
                0, & \text{otherwise}.
            \end{cases}
        \end{align*}
        Now, we consider the following differential equation:
        \begin{align}\label{eq:Proof1_Eq1}
            \frac{\partial B(\mathbf{x}(t))}{\partial t} = \lambda B(\mathbf{x}(t)) + h_p(\mathbf{x}(t)). 
        \end{align}        
        For an absolutely continuous solution $\vartheta_{\bar{x}}\colon [0,\infty) \rightarrow \R^n$ with initial value $\vartheta_{\bar{x}}(0) = \bar{x}$ it holds
        \begin{align}\label{eq:Proof1_Eq2}
            B(\vartheta_{\bar{x}}(t)) = \left(\int_0^t h_p(\vartheta_{\bar{x}}(\tau))e^{-\lambda \tau} d \tau + B(\bar{x}) \right) e^{\lambda t}.
        \end{align}
        Note that $B(\vartheta_{\bar{x}}(t)) \leq 0$ if $B(\bar{x})\leq0$ since $h_p(x) \leq 0$ for all $x \in X$ and for all $p \in \{1,\dots,N\}$ by \eqref{eq:Barrier_dot}.
        
        Now, we consider the intervals $[\tilde{t}_k,\tilde{t}_{k+1})$ for which $\boldsymbol{\sigma}(t)$ is constant and $(\tilde{t}_k,\tilde{t}_{k+1}) \cap \Omega = \emptyset$, i.e., the system dynamics of \eqref{eq:system_control} is just given by the differential equation ${\dot{x}(t) = f_p(\mathbf{x}(t),\mathbf{u}(t))}$, where $p=\boldsymbol{\sigma}(t)$ for all ${t \in [\tilde{t}_k,\tilde{t}_{k+1})}$ (the system does not undergo a switch or an impulse).
        Hence, on these intervals, a solution $\varphi_{x_0,\mathbf{u}, \boldsymbol{\sigma}}$ of \eqref{eq:system_control} coincides with a solution $\vartheta_{\bar{x}}$ of \eqref{eq:Proof1_Eq1} with initial value ${\vartheta_{\bar{x}}(0) = \bar{x} = \varphi_{x_0,\mathbf{u}, \boldsymbol{\sigma}}(\tilde{t}_k)}$ and hence, is given by \eqref{eq:Proof1_Eq2}.
        Thus, if $B(\bar{x}) = B(\varphi_{x_0,\mathbf{u}, \boldsymbol{\sigma}}(\tilde{t}_k)) \le 0$, also ${B(\varphi_{x_0,\mathbf{u}, \boldsymbol{\sigma}}(t)) \le 0}$ for all $t \in [\tilde{t}_k,\tilde{t}_{k+1})$.
        For ${\tilde{t}_0 = 0}$, this condition is satisfied by assumption, cf.\ \eqref{eq:Barrier_safe}. 
        If ${\boldsymbol{\sigma}(\tilde{t}_1) \neq \boldsymbol{\sigma}(\tilde{t}_0)}$, the system undergoes a switch in $\tilde{t}_1$. Thus, $\varphi_{x_0,\mathbf{u}, \boldsymbol{\sigma}}(\tilde{t}_1) = \lim_{t\nearrow \tilde{t}_1} \varphi_{x_0,\mathbf{u}, \boldsymbol{\sigma}}(t)$ and hence $B(\varphi_{x_0,\mathbf{u}, \boldsymbol{\sigma}}(\tilde{t}_1)) \le 0$. 
        If $\tilde{t}_1 \in \Omega$, the system undergoes an impulse in $\tilde{t}_1$. Consequently, 
        \[
        \varphi_{x_0,\mathbf{u}, \boldsymbol{\sigma}}(\tilde{t}_1) = 
        g\bigg(
        \underbrace{\lim_{t\nearrow \tilde{t}_1} \varphi_{x_0,\mathbf{u}, \boldsymbol{\sigma}}(t)}_{=:\tilde{x} \in \overline{\safe}}, 
        \underbrace{\lim_{t\nearrow \tilde{t}_1} \mathbf{u}(\varphi_{x_0,\mathbf{u}, \boldsymbol{\sigma}}(t))}_{=\mathbf{u}(\tilde{x})},
        \tilde{t}_1\bigg),
        \] 
        and therefore, $B(\varphi_{x_0,\mathbf{u}, \boldsymbol{\sigma}}(\tilde{t}_1)) \le \theta B(\tilde{x})$ by \eqref{eq:Barrier_impulse}.
        Since ${B(\tilde{x}) = B(\lim_{t\nearrow \tilde{t}_1} \varphi_{x_0,\mathbf{u}, \boldsymbol{\sigma}}(t))}$ and $\varphi_{x_0,\mathbf{u}, \boldsymbol{\sigma}}(t)$ for all ${t \in [0,\tilde{t}_1)}$, it holds $B(\tilde{x}) \le 0$, i.e., $B(\varphi_{x_0,\mathbf{u}, \boldsymbol{\sigma}}(\tilde{t}_1)) \le 0$.
        By an induction argument, we get $B(\varphi_{x_0,\mathbf{u}, \boldsymbol{\sigma}}(t))\le 0$ for all $t \in [0,\infty)$ and by \eqref{eq:Barrier_unsafe} the system is safe.        
    \end{pf}

    \begin{rem}
        Barrier functions are often only defined for ${\lambda \in (-\infty,0]}$. Allowing for positive values of $\lambda$ leads to even stricter conditions (as $B(x) \le 0$ for $x \in \safe$) and, thus, is not an issue. See also \cite{KHS+2013}. 
    \end{rem}

\subsection{Interconnected Control Systems}\label{subsec:InterconnectedSystems}
    In order to formally state what we mean by an interconnected control system, we first define the subsystems that comprise it. More specifically, we define the switched impulsive subsystems, based on \cite{alaoui2023symbolic}.
    \begin{defn}\label{def:subsystem}
         \emph{A switched impulsive subsystem} $\Sigma_i$, ${i \in \{1,\dots,M\}}, M \in \N$, is a tuple
        \begin{align}\label{eq:subsystem}
            \Sigma_i = (X_i, W_i, U_i,  \Omega_i, (f_{i,p})_{p \in \{1,\dots,N\}}, g_i)
        \end{align}
        with $X_i\subseteq \R^{n_i}$, $U_i \subseteq \R^{n^u_i}$, and $\Omega_i \subseteq \R_{\ge 0}$ as in Definition~\ref{def:system}.
        By $W_i \subseteq \R^{n^{\omega}_i}$, we denote the \emph{set of the internal inputs}.  
        Now, the family of transition maps ${f_{i,p}\colon X_i \times W_i \times U_i \rightarrow X_i}$, $p \in \{1,\dots, N\}$, and the jump map ${g_i\colon X_i \times W_i \times U_i \times \R_{\ge 0} \rightarrow X_i}$ also depend on the internal input $\omega_i \in W_i$. 
        Hence, the dynamics is given by 
        \begin{align}\label{eq:subsystem_description}
        \Sigma_i\colon\begin{cases}
          \dot{\mathbf{x}}_{i}(t) = f_{i,\boldsymbol{\sigma}(t)}(\mathbf{x}_i(t),\boldsymbol{\omega}_i(t),\mathbf{u}_i(t)), &\!\!\text{if } t \in  \R_{\geq 0}\!\setminus\! \Omega_i,\\
          \mathbf{x}_i(t) = g_i(\mathbf{x}_i^-(t),  \boldsymbol{\omega}_i^-(t), \mathbf{u}_i^-(t),t), &\!\!\text{if } t\in \Omega_i.\\
        \end{cases}
        \end{align}
        As before, $\mathbf{u}_i\colon [0,\infty)\rightarrow U_i$ and ${\boldsymbol{\sigma}\colon [0,\infty) \rightarrow \{1,\dots,N\}}$ are the control input and the admissible switching signal, respectively. 
        Furthermore, we assume that the internal input ${\boldsymbol{\omega}_i\colon [0,\infty)\rightarrow W_i}$ is absolutely continuous.
    \end{defn}

    Based on this definition, the interconnected switched impulsive system can be constructed. Note that the switching signal ${\boldsymbol{\sigma}}$ is not indexed by $i$, i.e., it is shared over all subsystems. This is an ease of notation but not a limitation, as discussed in Remark~\ref{rem:swichting_singal}.
    
    \begin{defn}\label{df:overallsystem}
        Let $\Sigma_i=(X_i, W_i, U_i, \Omega_i, (f_{i,p})_{p \in \{1,\dots,N\}}, g_i)$, $i \in \{1,\dots, M\}$, be $M$ switched impulsive subsystems as in Definition~\ref{def:subsystem} with 
        \[W_i = \prod_{\substack{j \in I_i}} X_{j}, \quad \text{for all}~i \in \{1,\dots,M\},\]
        where $I_i\subseteq \{1,\dots,M\}$ with $i \notin I_i$, i.e., elements $\omega_i \in W_i$ can be written as 
        $\omega_i=[x_{i_1},\dots,x_{i_{M_i}}]$
        with $M_i$ being the number of elements in $I_i$ and $x_{i_j} \in X_{i_j}$ for $i_j \in I_i$, ${j \in \{1,\dots, M_i\}}$. 
        The \emph{interconnected switched impulsive system} ${\Sigma=\mathcal{I}(\Sigma_1,\dots,\Sigma_M)}$ consisting of these $M$ subsystems $\Sigma_i$ is a switched impulsive system, i.e., a tuple
        \begin{align}\label{eq:overallsystem}
            \Sigma=(X,U,\Omega, (f_p)_{p\in \{1,\dots,N\}}, g),
        \end{align}
        where 
        \begin{align*}
            f_p(x,u)&=(f_{1,p}(x_1,\omega_1,u_1), \ldots, f_{M,p}(x_M,\omega_M,u_M))^{\top}\!, \\
            g(x,u,t)&=(\beta_1(x_1,\omega_1,u_1,t), \ldots, \beta_M(x_M,\omega_M,u_M,t))^{\top},
        \end{align*}
        with
        \begin{align*}
            \beta_i(x_i,\omega_i,u_i,t) =
            \begin{cases}
                x_i, \quad &\text{if } t \in \R_{\ge 0}\setminus\Omega_i,\\
                g_i(x_i, \omega_i,u_i,t), \quad &\text{if } t \in \Omega_i,
            \end{cases}
        \end{align*}
        and $X\subseteq\prod_{i=1}^{M}X_i$, $U\subseteq\prod_{i=1}^{M}U_i$ and $\Omega=\bigcup_{i=1}^{M}\Omega_i$.
    \end{defn}

    \begin{rem} \label{rem:swichting_singal}
        As stated before, all subsystems $\Sigma_i$ share the same (one-dimensional) switching signal ${\boldsymbol{\sigma}}$, i.e.,  the points in time at which the system dynamics can potentially change are the same for all subsystems. 
        This is no restriction. Assume we have a system consisting of two subsystems, $\Sigma_1$ and $\Sigma_2$. The first subsystem changes the dynamics at $t_1$ and the second one at $t_2 > t_1$. Then we can define three different modes, $p_1,p_2,p_3$, for the global system, where $f_{2,p_1} = f_{2,p_2}$ and $f_{1,p_2} = f_{2,p_3}$, i.e., for the global modes $p_1$ and $p_2$ the system dynamics do not change in $\Sigma_2$ and the same for $\Sigma_1$ and $p_2$ and $p_3$.
        The shared switching signal is given by $\boldsymbol{\sigma}(t) = p_1$ if $t \in [0,t_1)$, $\boldsymbol{\sigma}(t) = p_2$ if $t \in [t_1,t_2)$ and $\boldsymbol{\sigma}(t) = p_3$ if $t \in [t_2,\infty)$.
    \end{rem}

%% file: results.tex
\section{Safety condition for the interconnected system} \label{sec:Main}

    In this section, we will establish sufficient conditions to verify the safety of an interconnected and controlled switched impulsive system using barrier functions. Since we assume that the system is a large interconnected system, our goal is to develop a method that allows the compositional construction of barrier functions.
     
    From now on, let $\Sigma=\mathcal{I}(\Sigma_1,\dots,\Sigma_M)$ be an interconnected switched impulsive system consisting of $M$ switched impulsive subsystems $\Sigma_i$.
    To derive a method for the compositional construction of barrier functions, we need to make some assumptions for the smaller subsystems. More precisely, we assume the existence of so-called \emph{local pseudo barrier functions} $B_i$ for the single subsystems $\Sigma_i$, as defined in Definition~\ref{def:localBarrier}. The global barrier function $B$ for the interconnected system $\Sigma$ can then be derived by a weighted-sum or max formulation as stated in Theorem~\ref{th:globalBarrier}. 
   
    \begin{defn}\label{def:localBarrier}
        A $C^1$-function $B_i\colon X_i \to \R$ is called a \emph{local pseudo barrier function} for the {subsystem~$\Sigma_i$} with $X_{0,i}, X_{\mathsf{unsafe},i}\subseteq X_i$ if a control ${\mathbf{u_i}\colon X_i \to U_i}$ exists with
        \vspace*{-1em}
        \begin{subequations}
        \label{eq:localBarrier}
        \begin{align}
            &B_i(x_i) \le \varepsilon_{1,i}, \quad\forall x_i \in X_{0,i}     \label{eq:localBarrier_safe}\\         
            &B_i(x_i) > \varepsilon_{2,i}, \quad\forall x_i \in X_{\mathsf{unsafe},i} \label{eq:localBarrier_unsafe}\\        
            &\frac{\partial B_i}{\partial x_i} f_{i,p}(x_i, \omega_i, \mathbf{u}_i(x_i)) \leq \lambda_i B_i(x_i) + \sum_{j\in I_i}\gamma_{ij} B_j(x_j),  \nonumber\\ 
                &\quad\quad\forall x \in X, \forall p \in \{1,\dots,N\}, \label{eq:localBarrier_dot}\\
            &B_i(g_i(x_i,\omega_i,\mathbf{u}_i(x_i),t)) \le  B_i(x_i), \forall \nonumber t \in \Omega_i~\text{and}\\
            &\qquad \forall x_i=P_i(x), w_i = (P_j(x))_{j \in I_i}~\text{with}~x \in \overline{\safe}, \label{eq:localBarrier_impulse}
        \end{align}
        \end{subequations}
        for some $\lambda_i,\varepsilon_{1,i},\varepsilon_{2,i} \in \R$ and $\gamma_{ij} \in \R^{\ge0}$, $i \in \{1,\dots,M\}$ and $j \in I_i$. Here, $P_i\colon X \to X_i$ denotes the projection to $X_i$, i.e., $P_i([x_1,\dots,x_M]) = x_i$.
    \end{defn}

    Note that \eqref{eq:localBarrier_dot} encompasses the local barrier functions from all other connected subsystems. Hence, an independent construction of the local barrier functions is not possible. However, we will first show that these conditions are sufficient to construct a global barrier function in the following theorem and give a more restrictive alternative that allows for independent construction of local barrier functions in Theorem~\ref{th:Adapted_Inequality}.
        
    \begin{thm}\label{th:globalBarrier}
        Let $\mathbf{u}\colon X \to U$ be a control law such that for each subsystem $\Sigma_i$ the conditions of the local pseudo barrier function $B_i\colon X_i \to \R$ in Definition~\ref{def:localBarrier} are satisfied for $\mathbf{u}_i$.
        Define the matrices
        \begin{equation}
            \Lambda := \operatorname{diag}(\lambda_1,\dots,\lambda_M) \quad \text{and} \quad \Gamma := (\gamma_{ij})_{i,j=1,\dots,M} 
        \end{equation}
        where $\gamma_{ij} = 0$ if $j \notin I_i$, $i \in \{1,\dots,M\}$, and
        \begin{equation}\label{eq:eta}
           \eta = \min_{i\in \{1,\dots,M\}} \lambda_i.
        \end{equation}
        Moreover, assume that $(\Lambda + \Gamma - \eta I)$ is irreducible. \\
        Now, let $k = (k_1, \cdots, k_M)^{\top} \in (\R^{> 0})^M$ be a left eigenvector and $v = (v_1, \cdots, v_M)^{\top} \in (\R^{> 0})^M$ a right eigenvector of $(\Lambda + \Gamma - \eta I)$, i.e., 
        \begin{align*}
            k^{\top}(\Lambda + \Gamma - \eta I) &= \nu k^{\top} \quad \text{and} \quad
            (\Lambda + \Gamma - \eta I)v &= \mu v
        \end{align*}
        for some $\nu,\mu \in \R^{>0}$, where $I$ denotes the identity matrix.\\
        Consider the following two cases:
        \begin{enumerate}[left=0em]
            \item  Assume that $X_{0,i} = P_i(X_0)$, $X_{\mathsf{unsafe,i}}=P_i(X_{\mathsf{unsafe}})$ for all $i \in \{1,\dots,M\}$ and              \begin{equation}\label{eq:sumK}
                    \sum_{i=1}^M k_i\varepsilon_{1,i} \le 0 \quad \text{and} \quad
                     \sum_{i=1}^M k_i\varepsilon_{2,i} \ge 0.
                \end{equation}
                Then, there exists a barrier function $\BsumP$ for the interconnected switched impulsive systems $\Sigma$ as in Theorem~\ref{th:barrier_basic} that can be computed by
                \begin{align}\label{eq:Bsum}
                    \Bsum{x} =  \sum_{i=1}^{M} k_i B_i(x_i) = k^{\top} B_{\mathsf{vec}}(x), 
                \end{align}
                where $B_{\mathsf{vec}}(x) = (B_1(x_1),\dots, B_M(x_M))^{\top}$.
            \item  
            Assume that  
            $\displaystyle X_0 \subseteq \prod_{i=1}^{M} X_{0,i}$,
            \[\unsafe \subseteq \bigcup_{i=1}^M X_1 \times \!\cdots \!\times X_{i-1} \times X_{\mathsf{unsafe},i} \times X_{i+1} \times \!\cdots \! \times X_M\]
            and $\varepsilon_{1,i}=\varepsilon_{2,i}=0$  for all $i \in \{1,\dots,M\}$.
    
            Then, there exists a barrier function $\BmaxP$ for the interconnected switched impulsive systems $\Sigma$ as in Theorem~\ref{th:barrier_basic} that can be computed by
            \begin{align}\label{eq:Bmax}
                \Bmax{x} =  \max_{i \in \{1,\dots,M\}} \frac{1}{v_i} B_i(x_i). 
            \end{align}
        \end{enumerate}

        Hence, according to Theorem~\ref{th:barrier_basic}, the system $\Sigma$ is safe with respect to the control function $\mathbf{u}\colon X\to U$ and the switching signal $\boldsymbol{\sigma}\colon [0,\infty)\to \{1,\dots,N\}$ in both cases.
        
    \end{thm}

    \begin{pf}
        First, note that due to the definition of $\eta$, cf.\ \eqref{eq:eta}, the matrix $A :=(\Lambda + \Gamma - \eta I)$ has only non-negative entries and is irreducible by assumption. Thus, according to the Perron-Frobenius Theorem \cite[Chapter XIII.3, Theorem 2]{gantmacher1960a}, there exists a left eigenvector $k^{\top}$ and a right eigenvector $v$ that have only positive entries, i.e., $k_i,v_i > 0$, $i \in \{1,\dots,M\}$. 
        
        To prove the theorem, the four inequalities \eqref{eq:Barrier_safe}--\eqref{eq:Barrier_impulse} have to be verified for $\BsumP$ and $\BmaxP$. 
        It is easy to see that the inequalities \eqref{eq:Barrier_safe} and \eqref{eq:Barrier_unsafe} clearly hold in both cases. 
        Furthermore, in the first case, we have for all $x \in \overline{\safe}$:
        \begin{align*}
           &\Bsum{g(x,\mathbf{u}(x),t)}
            = \sum_{i=1}^{M} k_i B_i(\beta_i(x_i,\omega_i,\mathbf{u}_i(x_i),t)) \\
            &=\!\!\sum_{i:t\in \Omega_i} k_i B_i(g_i(x_i,\omega_i,\mathbf{u}_i(x_i),t))\!+\! \sum_{i:t\notin \Omega_i} k_i B_i(x_i) \\
            &\le\!\!\sum_{i:t\in \Omega_i} k_i B_i(x_i)\!+\!\!\!\sum_{i:t\notin \Omega_i}\!\!k_i B_i(x_i)
            = \Bsum{x}.
        \end{align*}
        Analogously, for the second case, we can show that $\Bmax{g(x,\mathbf{u}(x),t)} \le \Bmax{x}$ for all $x \in \overline{\safe}$, i.e., \eqref{eq:Barrier_impulse} is satisfied (with $\theta= 1$) as well.
 
        Finally, we prove \eqref{eq:Barrier_dot}. We start with the first case. For almost all $x \in X$ and all $p \in \{1,\dots,N\}$, we have:
        \begin{align*}
            \frac{\partial B}{\partial x}f_p(x,\mathbf{u}(x))  
            &= \sum_{i=1}^M k_i \frac{\partial B_i}{\partial x_i}(x_i)  f_{i,p}(x_i,\omega_i,\mathbf{u}_i(x_i))\\
            &\!\! \overset{\eqref{eq:localBarrier_dot}}{\le} \sum_{i=1}^M k_i \bigg(\lambda_i B_i(x_i) + \sum_{j \in I_i} \gamma_{ij}B_j(x_j)\bigg)\\
            &\!= k^{\top}\!(\Lambda + \Gamma) B_{\mathsf{vec}}(x) \\
            &\!= k^{\top}\!(\Lambda + \Gamma - \eta I) B_{\mathsf{vec}}(x) + \eta k^{\top}\!B_{\mathsf{vec}}(x)\\
            &\!= (\nu + \eta) B(x).
        \end{align*}
        Thus, \eqref{eq:Barrier_dot} holds in the first case.

        Since the maximum of $C^1$-functions is $C^1$ almost everywhere \citep{S2012}, we only have to consider the single functions $v_i^{-1}B_i(x_i)$ to prove \eqref{eq:Barrier_dot} in the second case. Thus, for $x \in X$, we assume that $\displaystyle\max_{i \in \{1,\dots,M\}} v_i^{-1}B_i(x_i)$ is uniquely determined by $v_l^{-1} B_l(x_l)$ for some $l \in \{1,\dots,M\}$. 
        Then, for all $p \in \{1,\dots,N\}$, we have:
        \begin{align*}
            \frac{\partial B}{\partial x} &f_p(x,\mathbf{u}(x)) 
            \!=  \frac{1}{v_l} \frac{\partial B_l}{\partial x_l}(x_l)  f_{l,p}(x_l,\omega_l,\mathbf{u}_l(x_l))\\
            &\!\! \overset{\eqref{eq:localBarrier_dot}}{\le} \frac{1}{v_l} \Big(\lambda_l B_l(x_l) + \sum_{j\in I_l}\gamma_{lj} B_j(x_j)\Big) \\
            &\!= \frac{1}{v_l} \Big( \lambda_l B_l(x_l) + \sum_{j\in I_l}\underbrace{\gamma_{lj} v_j \vphantom{\frac{1}{v_j}}}_{\ge 0} \underbrace{\frac{1}{v_j}B_j(x_j)}_{\le \frac{1}{v_l} B_l(x_l)}\Big)\\
            &\!\le \frac{1}{v_l} \Big(\lambda_l  B_l(x_l) + \frac{1}{v_l} B_l(x_l) \sum_{j\in I_l} \gamma_{lj} v_j \Big)\\
            &\!= \frac{1}{v_l}\Big(v_l\lambda_l +\sum_{j\in I_l} \gamma_{lj} v_j\Big)B(x) \\
            &\!= \frac{1}{v_l} \left[\left(\Lambda + \Gamma -\eta I\right) v + \eta v \right]_l B(x)\\
            &\!= \frac{1}{v_l} \left(\mu v_l + \eta v_l \right) B(x) = (\mu + \eta) B(x).
        \end{align*}
    \end{pf}

\begin{rem}
    Condition \eqref{eq:sumK} is hard to verify since it can only be checked after the overall construction. 
    However, we need to allow for this relaxation in \eqref{eq:localBarrier_safe} and \eqref{eq:localBarrier_unsafe}, i.e., $\varepsilon_{1,i}$ may be greater $0$ and $\varepsilon_{2,i}$ smaller $0$, since the projection of $\safe$ and $\unsafe$ to $X_i$ might intersect. Fortunately, in many practical applications, this is not the case, and it holds $X_0 = P_1(X_0) \times \dots \times P_M(X_0)$, $\unsafe = P_1(\unsafe) \times \dots \times P_M(\unsafe)$ and $P_i(X_0) \cap P_i(\unsafe) = \emptyset$. In this case, we can choose $\varepsilon_{1,i} = \varepsilon_{2,i}= 0$ and \eqref{eq:sumK} follows directly.
\end{rem}

\begin{figure}[t]
    \centering
    \begin{subfigure}[b]{0.36\textwidth}
        \begin{overpic}[width=\textwidth]{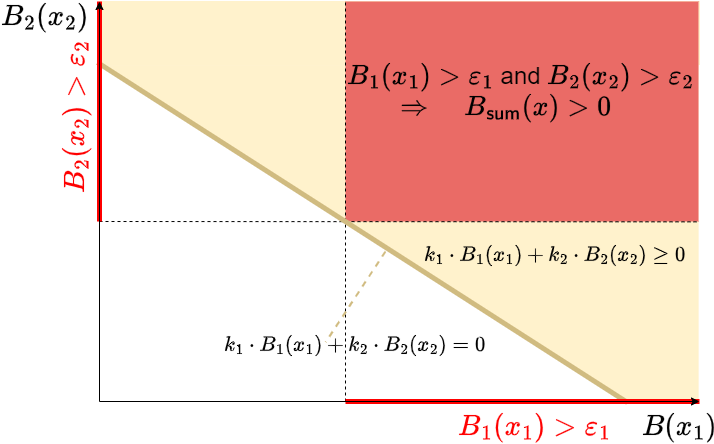}
            \put(-10,58){(a)}
        \end{overpic}
    \end{subfigure}
    \ \\[.5em]
    \begin{subfigure}[b]{0.36\textwidth}
        \begin{overpic}[width=\textwidth]{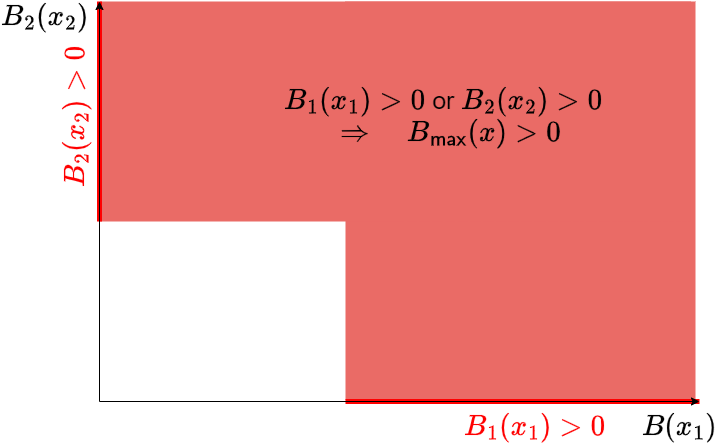}
            \put(-10,58){(b)}
        \end{overpic}
    \end{subfigure}
    \caption{Illustration of Case 1) and 2) of Theorem~\ref{th:globalBarrier} in (a) and (b), respectively. (b) For $\BmaxP$, it can be ensured that the complete red area is not entered since $\Bmax{x} > 0$. (a) For $\BsumP$, the same can be ensured for the whole beige and red area.
    However, the beige area is determined by $k$, and thus, only the red area can be determined a priori.
    }\label{fig:Theorem_cases}
    \vspace*{-1.5\baselineskip}
\end{figure}

\begin{rem} \label{rem:diff-sum-max}
    The most important difference between Case 1) and 2) in Theorem~\ref{th:globalBarrier} is the choice of the unsafe set $\unsafe$. For the max-formulation, i.e., for Case 2), we need to assume%
    \begin{equation*}
    \unsafe \subseteq \bigcup_{i=1}^M X_1 \times \dots \times X_{i-1} \times X_{\mathsf{unsafe},i} \times X_{i+1} \times \dots X_M.
    \end{equation*}
    This essentially means that if one subsystem $\Sigma_i$ is classified as unsafe, the global system $\Sigma$ is already unsafe. 
    On the contrary, for such a system, it holds that $P_i(X_{\mathsf{unsafe}}) = \R$, making the sum-formulation infeasible due to \eqref{eq:localBarrier_safe} and \eqref{eq:localBarrier_unsafe}. However, in the case that the global system $\Sigma$ is only unsafe if all the subsystems $\Sigma_i$ are classified as unsafe, i.e., 
    \[\unsafe \subseteq X_{\mathsf{unsafe},1} \times \dots \times X_{\mathsf{unsafe},M},\] the sum-formulation is beneficial, as $P_i(X_{\mathsf{unsafe}})=X_{\mathsf{unsafe},i}$, and the max-formulation is not applicable. 
    This is illustrated in Fig.~\ref{fig:Theorem_cases} and discussed for the numerical example in Section~\ref{sec:Numerics}.
\end{rem}

Since condition \eqref{eq:localBarrier_dot} includes all local barrier functions of connected subsystems, the local barrier functions cannot be constructed independently. Thus, we provide alternative conditions, which, however, limit the systems to which our results are applicable.

\begin{thm}\label{th:Adapted_Inequality}
    Assume that the conditions of Theorem~\ref{th:globalBarrier} are satisfied but replace \eqref{eq:localBarrier_dot} in Definition~\ref{def:localBarrier} by the following inequalities that have to hold for each subsystem $\Sigma_i$, ${i \in \{1,\dots,M\}}$:
        \begin{subequations}
            \label{eq:localBarrier_adapt}
            \begin{align}
                & \|x_i\|_{2}^2 \le \alpha_i B_i(x_i) + L, \forall x_i \in X_i \nonumber\\
                &\quad \text{if}~i \in I_j~\text{ for some}~j \in \{1,\dots,M\}~\text{with}~c_j \neq 0 \label{eq:localBarrier_new}, \\ %
                & \frac{\partial B_i}{\partial x_i} f_{i,p}(x_i, \omega_i, \mathbf{u}_i(x_i)) \leq \lambda_i B_i(x_i) + c_i\left( \|{\omega_i}\|_2^2 - LM_i \right), \nonumber\\%
                & \quad \forall x_i \in X_i, \omega_i \in W_i, \forall p \in \{1,\dots,N\}, 
                \label{eq:localBarrier_dot_adapt}
            \end{align}
        \end{subequations}        
    where $c_i = 0$ if $M_i= 0$, $M_i$ denotes the number of indices in $I_i$ and $c_i, \alpha_i \in \R^{\ge0}$ and $L > 0$ are some constants, whereby the latter has to be the same for all subsystems. Define $\gamma_{ij}:=c_i\alpha_j$ and let $k$ be the positive eigenvector as before.
    Then, \eqref{eq:localBarrier_dot} is satisfied and, thus, Theorem~\ref{th:globalBarrier} holds, i.e., the system $\Sigma$ is safe. 
\end{thm}

\begin{pf}
    For all $x \in X$ and $\forall p \in \{1,\dots,N\}$, we have:
        \begin{align*}
            &\frac{\partial B}{\partial x}f_p(x,\mathbf{u}(x))
            = \sum_{i=1}^M k_i \frac{\partial B_i}{\partial x_i}(x_i)  f_{i,p}(x_i,\omega_i,\mathbf{u}_i(x_i))\\
            &\!\!\overset{\eqref{eq:localBarrier_dot_adapt}}{\le} \sum_{i=1}^M k_i \left(\lambda_i B_i(x_i) + c_i\left( \|{\omega_i}\|_2^2 - LM_i \right)\right)\\
            &= \sum_{i=1}^M k_i \bigg(\lambda_i B_i(x_i) + c_i \Big( \sum_{j \in I_i} \|x_j\|_2^2 - LM_i\Big)\bigg)\\
            &= \sum_{i=1}^M k_i \bigg(\lambda_i B_i(x_i) + \sum_{j \in I_i} c_i \left( \|x_j\|_2^2 - L\right)\bigg)\\
            &\!\!\overset{\eqref{eq:localBarrier_new}}{\le} \sum_{i=1}^M k_i \bigg(\lambda_i B_i(x_i) + \sum_{j \in I_i} \gamma_{ij} B_j(x_j) \bigg),
        \end{align*}
    where we used in the last inequality that $k_i, c_i \ge 0$
    and $\gamma_{ij} = c_i \alpha_j$ for $i \in \{1,\ldots,M\}$,  $j \in I_i$.
\end{pf}

%% file: simulation.tex
\section{Numerical Example}\label{sec:Numerics}

To illustrate the results of the paper, we consider an epidemiological model based on \cite{yang2010effects} and \cite{liu2017infectious}.
More precisely, we utilize a pulse vaccination model that takes into account the population movement across a network of $n$ patches and incorporates the influence of seasonal variations on disease transmission dynamics. 
Such a model can be described by \eqref{eq:subsystem_description} with $x_{i} \coloneqq (S_i, I_i, R_i)$, where $S_i$, $I_i$, and $R_i$ are the number of susceptible, infectious, and recovered individuals in the \mbox{$i$-th} patch, respectively, and $\Omega_i\coloneqq\Omega = \{kT : k \in \N\}$, where $T > 0$ represents the vaccination interval. 
The dynamics in the $i$-th patch is given by the subsystem %
\begin{equation*}
\label{eq:numerical_example}
\!\Sigma_i\!\!:\!\!
\begin{cases}
    {\text{For } t \in \mathbb{R}_{\ge 0} \setminus \Omega:} \\
    \dot{\mathbf{S}}_{i}(t)\!=\!m_{i}\!-\!\beta_{i,\boldsymbol{\sigma}(t)}\mathbf{S}_{i}(t)\mathbf{I}_{i}(t)\!-\!\mu_{i}\mathbf{S}_{i}(t)\!\\ \qquad \qquad +\!\textstyle\sum_{j=1}^{M}a_{ij}\mathbf{S}_{j}(t), \\ 
    \dot{\mathbf{I}}_{i}(t)\!=\!\beta_{i,\boldsymbol{\sigma}(t)}\mathbf{S}_{i}(t)\mathbf{I}_{i}(t)\!-\!(\mu_{i}+r_{i})\mathbf{I}_{i}(t)\!\\ \qquad \qquad +\!\textstyle\sum_{j=1}^{M}b_{ij}\mathbf{I}_{j}(t), \\
    \dot{\mathbf{R}}_{i}(t)\!=\!r_{i}\mathbf{I}_{i}(t)\!-\!\mu_{i}\mathbf{R}_{i}(t)\!+\!\textstyle\sum_{j=1}^{M}c_{ij}\mathbf{R}_{j}(t), \\[0.25em]
    {\text{For } t \in \Omega:} \\
    \mathbf{S}_{i}(t)\!=\!(1-p_{i})\mathbf{S}^{-}_{i}(t), \\
    \mathbf{I}_{i}(t)\!=\!\mathbf{I}^{-}_{i}(t), \\
    \mathbf{R}_{i}(t)\!=\!\mathbf{R}^{-}_{i}(t)+p_{i}\mathbf{S}^{-}_{i}(t),
\end{cases}
\end{equation*}
where $p_i$ is the proportion of successfully vaccinated persons, $m_i$ is the constant recruitment rate, 
$\mu_i$ is the death rate, and $r_i$ is the recovery rate of infectious individuals. 
For $i \ne j$, $a_{ij}$, $b_{ij}$, and $c_{ij}$ represent the immigration rates of susceptible, infectious, and recovered individuals from the $j$-th to the $i$-th patch. Moreover, $-a_{ii}$, $-b_{ii}$, and $-c_{ii}$ represent the emigration rates of susceptible, infectious,
and recovered individuals from the $i$-th patch, respectively. 
The transmission rate $\beta_{i,\boldsymbol{\sigma}(t)}$ varies depending on the specific patch and can additionally be influenced by time-dependent factors, such as the initiation of a rainy season \citep{liu2017infectious}. This temporal variation is represented by the switching signal~$\boldsymbol{\sigma}(t)$.

For the first experiment, we consider a system with ${M=3}$ patches, where ${X=\{x \in \R^{3M}\colon ||x||_2\leq 100\}}$ and ${X_i = \{ x_i=(S_i,I_i,R_i) \in \R^{3}\colon ||x_i||_2\leq 100 \}}$, $i \in \{1,2,3\}$. 
Moreover, we choose $m_i=0.02$, $r_i=0.05$, $p_i= 0.8$, and the migration rates
\begin{align*}
 [a_{ij}] = 10^{-4}&\left[\begin{smallmatrix}
   -7 & 14 & 5 \\
    6 & -21 & 19 \\
    1 & 7 & -24 \\
\end{smallmatrix}\right],
 [b_{ij}] = 10^{-4}\left[\begin{smallmatrix}
   -14 & 2 & 2 \\
    11 & -5 & 1 \\
    3 & 3 & -3 \\
\end{smallmatrix}\right],\\
 &[c_{ij}] = 10^{-4}\left[\begin{smallmatrix}
   -2 & 3 & 4 \\
    1 & -6 & 13 \\
    1 & 3 & -17 \\
\end{smallmatrix}\right].
\end{align*}

We aim to validate the system's safety under an (uncontrolled) arbitrary switching $\boldsymbol{\sigma}: [0,\infty) \to \{1,2\}$ such that   
$\beta_{i,1} = 0.001$ and $\beta_{i,2} = 0.005$ for $i \in \{1,\dots,3\}$. 
We assume that the system is unsafe if $I_i\geq 5$ for at least one $i \in \{1,2,3\}$. 
More precisely, for all subsystems, we define 
\begin{align*}
    X_{\mathsf{unsafe},i} &= \{ x_i\in\R^3 \colon q_{1,i}(x_i) \ge 0, q_{2,i}(x_i) \le 0\},\\ 
    X_{\mathsf{safe},i} &= \{ x_i\in\R^3 \colon q_{1,i}(x_i) \ge 0, q_{2,i}(x_i) \ge 0\},~\text{and} \\
    X_{0,i} &= \{ x_i\in\R^3 \colon q_{3,i}(x_i) \ge 0\}, 
\end{align*}
where ${q_{1,i},q_{2,i},q_{3,i}:X_i\rightarrow \R}$ are given by
\begin{align*}
    q_{1,i}(x_i) &= 100 - (S_i^2 + I_i^2 + R_i^2),\\ 
    q_{2,i}(x_i) &=  5 - I_i \quad \text{and}\quad q_{3,i}(x_i) = 1 - (S_i^2 + I_i^2 + R_i^2).
\end{align*}
As the initial set of the global system we choose ${X_0 = X_{0,1} \times X_{0,2} \times X_{0,3}}$ and the unsafe set is given by
$\unsafe = (X_{\mathsf{unsafe},1} \times X_2 \times X_3) \cup \dots \cup (X_1 \times X_2 \times X_{\mathsf{unsafe},3})$.
This structure fits the second case of Theorem~\ref{th:globalBarrier} as explained in Remark~\ref{rem:diff-sum-max}.

To computationally construct local barrier functions, we assume a polynomial structure for each $B_i$ and use the sum-of-squares (SOS) framework. The goal is to find a polynomial $B_i$, and sum-of-squares polynomials ${\rho_{j,i}:X_i\rightarrow \R}$, $j={1,\ldots,5}$, and ${\tau_{k,i}:X_i\rightarrow \R}$, ${k=1,\ldots,7}$, such that  
\begin{align*}
   \rho_{1,i}(x_i)\!&=\! -B_i(x_i)\!+\!\varepsilon_{1,i}\!-\!\tau_{1,i}(x_i)q_{3,i}(x_i), \\
   \rho_{2,i}(x_i)\!&=\!\! B_i(x_i)\!-\!\varepsilon_{2,i}\!-\!\tau_{2,i}(x_i)q_{1,i}(x_i)\!+\!\tau_{3,i}(x_i)q_{2,i}(x_i), \\
   \rho_{3,i}(x_i)\!&=\!\! B_i(x_i)\!\!-\!\!B_i^{+}(x_i)\!\!-\!\tau_{4,i}(x_i)q_{1,i}(x_i)\!\!+\!\tau_{5,i}(x_i)q_{2,i}(x_i), \\
   \rho_{4,i}(x_i)\!&=\! -\|x_i\|_{2}^2\!+\!\alpha_i B_i(x_i)\!+\!L\!-\!\tau_{6,i}(x_i)q_{1,i}(x_i), \\
   \rho_{5,i}(x_i)\!&=\! -\frac{\partial B_i}{\partial x_i} f_{i,p}(x_i, \omega_i, \mathbf{u}_i(x_i)) \!+\! \lambda_i B_i(x_i) \\
    & \qquad\! +c_i\left( \|{\omega_i}\|_2^2 - LM_i \right) -\tau_{7,i}(x_i),
\end{align*}
with $B_i^{+}(x_i):=B_i(g_i(x_i))$, transforming the problem in a convex optimization problem \cite{parrilo2003semidefinite}. If such polynomials exist,
the conditions of Theorem~\ref{th:Adapted_Inequality} are satisfied, and the system is safe. Using the sum-of-squares parser of YALMIP \citep{Lofberg2009}, and choosing $L=10000$, $\lambda_i = -0.1$, $\alpha_i = 0.1$ and $c_i = 0.0002$, we can compute a (global) barrier function as in \eqref{eq:Bmax}, where 
\begin{align}\label{eq:barrier_for_example1}
\frac{B_1}{v_1}\!=\!&-3332.54\!-\!0.35 S_1\!+\!345.95I_1\!-\!0.35 R_1\!-\!38.58S_1 I_1 \nonumber\\ &\!\!\!\!-3.50 S_1^2\!+\!190.42 I_1^2\!-\!4.24 R_1^2\!-\!8.48 S_1 R_1\!-\!38.58I_1 R_1, \nonumber \\
\frac{B_2}{v_2}\!=\!&-3319.49\!+\!7.44 S_2\!+\!325.58 I_2\!+\!7.44 R_2\!-\!34.59 S_2 I_2 \nonumber\\ &\!\!\!\!-0.30 S_2^2\!+\!160.27 I_2^2\!-\!1.03 R_2^2\!-\!2.07 S_2 R_2\!-\!34.59 I_2 R_2, \nonumber \\ 
\frac{B_3}{v_3}\!=\!&-3313.49\!+\!6.43 S_3\!+\!334.78 I_3\!+\!6.43 R_3\!-\!35.77 S_3 I_3 \nonumber\\ &\!\!\!\!-0.38 S_3^2\!+\!162.16 I_3^2\!-\!1.11 R_3^2\!-\!2.23 S_3 R_3\!-\!35.77 I_3 R_3. \nonumber
\end{align}
with the right eigenvector $v=(0.57, 0.58, 0.58)^\top$, proving that the system is safe.\footnote{%
Verification of the validity of the SOS conditions are done by checking the residual of the primal problem (the computed residuals were less than $10^{-5}$) and using Theorem~4 of \cite{Lofberg2009}.
}
Note that in this setting, it is impossible to use Theorem~\ref{th:barrier_basic} with SOS for the computation of a global barrier function since the corresponding unsafe set is not a semi-algebraic set, showing that Theorem~\ref{th:globalBarrier} is more flexible regarding the definition of the safe and unsafe sets.

In the next example, we aim to demonstrate the performance increase in terms of the computation time that is possible by calculating the local barrier functions using Theorem~\ref{th:globalBarrier}, or more precisely Theorem~\ref{th:Adapted_Inequality}, instead of the global construction according to Theorem~\ref{th:barrier_basic}. 
Since we cannot compute the global barrier function using SOS in the previously described setting, we redefine the system to be unsafe if and only if the number of infected people exceeds the threshold in all subsystems, i.e., ${\unsafe=\prod_{i=1}^M X_{\mathsf{unsafe},i}}$, where $X_{\mathsf{unsafe},i}$, $i \in \{1,\dots,M\}$, is defined as before. 
Thus, for the construction of the global barrier function from the local barrier functions, we have to consider the sum formulation given by the first case of Theorem~\ref{th:globalBarrier}, cf.\ Remark~\ref{rem:diff-sum-max}. 
However, in this setting, $X_\mathsf{safe}$ is not an intersection of polynomial sub-level sets anymore. Fortunately, $X_\mathsf{safe}$ is only needed to verify \eqref{eq:Barrier_impulse}, which might be possible to check on the whole set~$X$. %

Despite the adaption of the safe and unsafe set, we consider the same parameters as before but with a flexible number of subsystems to study the computational time with increasing dimensions and complexity of the system. To this end, we assume that the patches are arranged in a ring structure, i.e., if $|i-j| \mod M>1$, then $a_{ij}$, $b_{ij}$ and $c_{ij}$ are set to $0$, and if $|i-j| \mod M=1$ holds, then $a_{ij}$, $b_{ij}$, and $c_{ij}$ are taken as the absolute value of a randomly sampled value from a normal distribution with $0$ mean and $0.001$ standard deviation.
To make sure that as many people emigrate as immigrate, we set $a_{ii}=-\sum_{i=1}^Ma_{ij}$ and similarly $b_{ii}=-\sum_{i=1}^Mb_{ij}$ and $c_{ii}=-\sum_{i=1}^Mc_{ij}$.

Using SOS, we search again for polynomial barrier functions of degree 2 satisfying the conditions of Theorem~\ref{th:barrier_basic} as well as for polynomial local pseudo barrier functions satisfying the conditions of Theorem~\ref{th:globalBarrier}. 
We run the experiment for different numbers of subsystems ranging from $M=3$ to $M=6$.
Our results are reported in Table~\ref{tab:comparison}, showing the mean and standard deviation of the running times calculated over five successful searches.

\begin{table}[ht]
    \caption[Comparison of computation times for global and local barrier functions]{Comparison of computation times for global and local barrier functions\footnotemark}\label{tab:comparison}
    \resizebox{\columnwidth}{!}{%
        \begin{tabular}{lcrr}
        \hline
            \multicolumn{1}{c}{$M$} &
            \multicolumn{1}{c}{\begin{tabular}[c]{@{}c@{}} 
                Dimension of \\ state space of $\Sigma$ 
            \end{tabular}} &
            \multicolumn{1}{c}{\begin{tabular}[c]{@{}c@{}}Run time \\ (Theorem~\ref{th:barrier_basic})\end{tabular}} &
            \multicolumn{1}{c}{\textbf{\begin{tabular}[c]{@{}c@{}}Run time \\ (Theorem~\ref{th:globalBarrier})\end{tabular}}} \\ 
        \hline
            3 & 9  & \textbf{1.67s $\boldsymbol{\pm}$ 0.18s}   & 2.10s $\pm$ 0.43s \\
            4 & 12 & 5.78s $\pm$ 0.16s   & \textbf{2.57s $\boldsymbol{\pm}$ 0.04s} \\
            5 & 15 & 24.66s  $\pm$ 0.18s & \textbf{3.19s $\boldsymbol{\pm}$ 0.10s} \\
            6 & 18 & 83.91s $\pm$ 3.83s  & \textbf{3.94s $\boldsymbol{\pm}$ 0.25s} \\ 
        \hline
        \end{tabular}%
    }
\end{table}

As expected, the running times increase along with the number of subsystems. However, the time for computing the global barrier function (Theorem~\ref{th:barrier_basic}) increases much faster than for the composed construction (Theorem~\ref{th:globalBarrier}).
In particular, compared to  Theorem~\ref{th:barrier_basic}, the verification of Theorem~\ref{th:globalBarrier} is almost twice as fast when dealing with 4 subsystems and about 20 times faster when dealing with~6. %